\documentclass[11pt]{article}
\usepackage{graphicx}
\usepackage{color}
\usepackage{amsmath,amsfonts,amssymb,amsthm}
\def\ds{\displaystyle}
\def\forall{\hbox{for all}~}
\def\L{{\bf L}}

\def\N{{\cal N}}

\def\ve{\varepsilon}
\def\n{\noindent}

\def\R{{\mathbb R}}

\def\tv{\mathop{\mathrm{Tot.Var.}}}
\def\vs{\vskip 2em}

\def\v{\vskip 1em}
\def\D{{\cal D}}
\def\O{{\cal O}}

\def\begi{\begin{itemize}}
\def\endi{\end{itemize}}
\def\C{{\cal C}}
\def\dint{\int\!\!\int}

\def\bega{\begin{array}}
\def\enda{\end{array}}

\def\bel{\begin{equation}\label}
\def\eeq{\end{equation}}
\def\sqr#1#2{\vbox{\hrule height .#2pt
\hbox{\vrule width .#2pt height #1pt \kern #1pt
\vrule width .#2pt}\hrule height .#2pt }}
\def\square{\sqr74}
\def\endproof{\hphantom{MM}\hfill\llap{$\square$}\goodbreak}

\newtheorem{theorem}{Theorem}[section]

\begin{document}
\title{\bf  A Remark on the Uniqueness of  Solutions to Hyperbolic Conservation Laws}
\vs

\author{Alberto Bressan$^*$ and  Camillo De Lellis$^{**}$
\\
\, \\
{\small (*) Department of Mathematics, Penn State University, University Park, Pa.16803, USA.}\,\\ 
{\small (**) Institute for Advanced Study, Princeton, NJ 08540, USA. } \\ \, \\
{\small E-mails:  axb62@psu.edu, ~ camillo.delellis@ias.edu}
\,\\}
\maketitle

\begin{abstract} Given a strictly hyperbolic $n\times n$ system of conservation laws, 
it is well known that there exists a unique Lipschitz semigroup of 
weak solutions, defined on a domain of functions with small total variation, which are limits
of vanishing viscosity approximations.
Aim of this note is to  prove that every weak solution taking values in the domain of the semigroup,
and whose shocks satisfy the  Liu admissibility conditions, actually 
coincides with a semigroup trajectory.
\end{abstract}
 \v
\section{Introduction}
\label{sec:1}
\setcounter{equation}{0}
We consider the Cauchy problem for  a strictly hyperbolic $n\times n$ system of conservation laws in one space dimension
\bel{1} u_t + f(u)_x ~=~0,\eeq
\bel{2} u(0,x)~=~\bar u(x),\eeq
with $f\in \C^2 (\mathbb R^n, \mathbb R^n)$.
In this setting, it is well known that there exists a Lipschitz continuous semigroup $S:\D\times [0, +\infty[\, \mapsto \D$ of weak solutions, 
defined on a domain $\D\subset\L^1(\R;\,\R^n)$ of  functions
with sufficiently small total variation.  The trajectories of this 
semigroup are  the unique limits of  
vanishing viscosity approximations \cite{BiB}.  All of their shocks satisfy the Liu admissibility conditions \cite{Bi, L1, L2}.  
We recall that the semigroup is globally Lipschitz continuous w.r.t.~the $\L^1$ distance. Namely,
there exists a constant $L$ such that
\bel{Slip}
\bigl\| S_t\bar u - S_s\bar v\bigr\|_{\L^1}~\leq~L\Big( |t-s| + \|\bar u - \bar v\|_{\L^1}\Big)\qquad
\forall s,t\geq 0,~~\bar u,\bar v\in\D.\eeq

Given any weak solution $u=u(t,x)$ of (\ref{1})-(\ref{2}),  
various conditions have been derived in \cite{BG, BLF, BL} which guarantee
the identity  
\bel{uniq} u(t)~=~S_t\bar u\qquad\forall t\geq 0.\eeq
Since the semigroup $S$ is unique, the identity
(\ref{uniq}) yields  the uniqueness of  solutions to the Cauchy problem (\ref{1})-(\ref{2}).
In addition to the standard assumption that each characteristic field is either linearly degenerate or genuinely nonlinear, earlier results required some additional
regularity conditions, such as ``Tame Variation" or ``Tame Oscillation", controlling
the behavior of the solution near a point where the variation is small.

Aim of the present note is to show that uniqueness is guaranteed in a fully general setting: 
without any assumption about genuine nonlinearity, and without
any of the above regularity conditions. Moreover, no assumption is made about the existence
of a convex entropy. Our only requirement is that all points of approximate jump
satisfy
the Liu admissibility conditions.

As in \cite{BG, BLF, BL}, 
the proof relies on the elementary error estimate
\bel{errest}
\bigl\| u(t)- S_t\bar u\bigr\|_{\L^1}~\leq~L\cdot\int_0^t \liminf_{h\to 0+}\,
{\bigl\|u(\tau+h)- S_h u(\tau)\bigr\|_{\L^1}\over h}\, d\tau\,.\eeq
Indeed,
we will prove that the integrand is zero for a.e.~time $\tau\geq 0$. 
Following an argument introduced in \cite{BGlimm}, this is achieved by two estimates:
\begi
\item[(i)] In a neighborhood of a point $(\tau,y)$ where $u(\tau,\cdot)$ has a large jump,
the weak solution $u$ is compared with the solution to a Riemann problem.
\item[(ii)] In a region where the total variation is small,
the weak solution $u$ is compared with the solution to a linear system with constant
coefficients.\endi

To fix ideas, let 
\bel{Mdef}M~\doteq~\sup\Big\{ \tv\{\bar u\,; \R\}\,;~~\bar u \in \D\Big\}\eeq
be an upper bound for the total variation of all functions in the domain of the semigroup.
Notice that this implies
\begin{equation}\label{e:L-infty}
\|\bar u\|_{\L^\infty}\leq M \qquad \forall \bar u\in {\cal D}\, . 
\end{equation}
Moreover, for each BV function $\bar u\in \D$, we shall take its right-continuous representative, so that
$\bar u(x)= \lim_{y\to x+} u(y)$.

To state our result, we first describe the basic setting.
\v\begi
\item[{\bf (A1)}] {\bf (Conservation equations)}  {\it The function $u=u(t,x)$
is a weak solution of the Cauchy problem (\ref{1})-(\ref{2}) taking values within the domain of the semigroup.\\
More precisely, $u:[0,T]\mapsto\D$ is continuous w.r.t.~the $\L^1$ distance.
The identity $u(0,\cdot)=\bar u$ holds  in $\L^1$, and moreover
\bel{4}\dint \big(u\varphi_t+f(u)\varphi_x\big) ~dxdt=0\eeq
for every $\C^1$ function $\varphi$ with compact support contained
inside the open strip $\, ]0,T[\,\times\R$.}
\endi
To introduce the Liu admissibility condition on the shocks \cite{L1, L2},  we first recall that 
{\bf (A1)} implies that $u=u(t,x)$ is a function of bounded variation in time and space (cf. \cite[Section 5.1]{EG} for the definition).
Indeed, by~\cite[Theorem~4.3.1]{D2}, we have
the Lipschitz bound
 \bel{ulip}
 \bigl\| u(t_2,\cdot)-u(t_1,\cdot)\bigr\|_{\L^1(\R)}~\leq~C_M\,(t_2-t_1)
 \qquad\qquad \forall 0\leq t_1\leq t_2\,,\eeq
 for some constant $C_M>0$ depending only on  the flux
 $f$ and on the upper bound $M$ for the total variation.
%
%

By the structure theorem for BV functions of two variables
(see e.g. \cite[Section 5.9]{EG} or \cite{AFP}), there is a Borel subset $\mathbf{J}\subset [0,1]\times \mathbb R$ with the following three properties.
\begin{itemize}
\item[(i)] Every point $(\tau, \xi)\not\in \mathbf{J}$ is a point of approximate continuity.
\item[(ii)] $\mathbf{J}$ is countably $1$-rectifiable, i.e. it can be covered by countably many Lipschitz curves, possibly leaving out a subset of zero $\mathcal{H}^1$ measure ($\mathcal{H}^1$ denotes the Hausdorff $1$-dimensional measure, cf. \cite[Section 2.1]{EG}).
\item[(iii)] $\mathcal{H}^1$-almost every point $(\tau,\xi)\in \mathbf{J}$ is an approximate jump of the function $u$. More precisely
there exist states 
$u^-,u^+\in\R^n$ and a speed $\lambda\in\R$ such that, calling
\bel{5}U(t,x)~\doteq~
\left\{\bega{rl} u^-\qquad &\hbox{if}\qquad x<\lambda t,
\\[1mm] u^+\qquad &\hbox{if}\qquad x>\lambda t,\enda\right.
\eeq
there holds
\bel{6}\lim_{r\to 0+}~{1\over r^2}
\int_{-r}^r\int_{-r}^r\Big| u(\tau+t,~\xi+x)-U(t,x)\Big|~dx\, dt~=~0.\eeq
\end{itemize}
Defining the rescaled functions
\begin{equation}\label{e:rescaled}
u_r (t,x) \,\doteq\, u (\tau + r t, \xi+r x)\, ,
\end{equation}
by (iii) and  {\bf (A1)} it follows that $u_r$ converges  to $U$ in $\L^1_{loc}(\R^2)$. In particular the conservation equations (\ref{4}) must hold for the 
piecewise constant function $U$, and the triple $(u^+, u^-, \lambda)$ must therefore satisfy the Rankine-Hugoniot equations:
\bel{RH}
f(u^+) - f(u^-)~=~\lambda (u^+-u^-).\eeq

Now let a left state $u^-\in \R^n$ be given.
Since the system is strictly hyperbolic, there exist $n$ shock curves
$$s\mapsto S_i(s)(u^-),\qquad\qquad i=1,\ldots,n,$$
parameterizing the sets of right states $u^+$ connected to the left state $u^-$ by 
a shock of the $i$-th family \cite{Bbook, D2, HR}.
As in  (\ref{RH}), denote by  $\lambda=\lambda_i(s)$ the Rankine-Hugoniot speed
of a shock with left and right states $u^-$ and $u^+=S_i(s)(u^-)$, 

\begi
\item[{\bf (A2)}] {\bf (Liu admissibility condition)} 
{\it In the above setting, a shock with left and right states 
$u^-$ and $u^+ = S_i(\sigma)( u^-)$ is {\bf Liu-admissible} if
$\lambda_i(\sigma) \leq \lambda_i(s)$ for all $s\in [0,\sigma]$.}
\endi

Our result can be simply stated as:
\begin{theorem}\label{t:1} Let (\ref{1}) be a strictly hyperbolic $n\times n$ system.
Then every weak solution $u:[0,T]\mapsto\D$,  taking values within the 
domain of the semigroup and whose shocks satisfy the Liu admissibility condition, 
coincides with a semigroup trajectory.  \end{theorem}

Under the additional assumptions that each characteristic family is either 
linearly degenerate or genuinely nonlinear, and  that the $n\times n$ system (\ref{1}) admits a strictly 
convex entropy selecting the admissible shocks, this
uniqueness result was recently proved in \cite{BGu}.  
Restricted to a class of 
$2\times 2$ systems, an earlier
proof  can also be found in \cite{CKV}.

\v

\section{Proof of the theorem}  
\label{sec:2}
\setcounter{equation}{0}

 {\bf 1.} Let  $\mathbf{J}$ be the set introduced in the previous section and let $\mathbf{S}\subset \mathbf{J}$ be the subset of all points which are not approximate jumps. 
 Since $\mathcal{H}^1 (\mathbf{S}) =0$, its projection on the time axis is a subset ${\cal N}\subset[0,T]$ which is null for the Lebesgue measure.
Every point $(t,x)\in [0,T]\times \R$ with $t\notin {\cal N}$ is therefore either a point
of approximate jump, or a point of approximate continuity.

Let us denote by ${\cal J}$ the set 
\[
{\cal J} \,\doteq\, \bigl\{(\tau, \xi): \tau\not \in {\cal N} \; \mbox{and $(\tau, \xi)$ is an approximate jump}\bigr\}\, .
\]
While it follows immediately from the aforementioned BV structure theorem that ${\cal J}$ is rectifiable, we claim here a stronger property: ${\cal J}$ can be covered 
by the graphs of countably many Lipschitz functions
\bel{p1} x=\phi_\ell(t),\qquad\qquad \ell\in {\mathbb N}\, ,\eeq
and moreover the Lipschitz constant of each $\phi_\ell$ is bounded by a number $\Lambda$ which depends only on $f$ and on the constant $M$ in (\ref{Mdef}). More precisely, by 
 recalling \eqref{e:L-infty} we can set 
 \bel{Ladef}\Lambda ~\doteq ~2 \,{\rm Lip} \bigl(f, B_M\bigr),\eeq
 where the right hand side denotes the Lipschitz constant of the function $f$ 
 over the ball $B_M\subset\R^n$ centered at the origin with  radius $M$. 

Since by (\ref{e:L-infty}) it follows
\[
|u^\pm (\tau, \xi)| \leq M \qquad \forall (\tau, \xi)\in {\cal J} ,
\]
the above definition implies that the shock speed $\lambda=\lambda (\tau, \xi)$ at (\ref{5}),
(\ref{RH})
satisfies the bound 
\begin{equation}\label{mspeed}
|\lambda (\tau, \xi)| ~\leq~ {\rm Lip}\, (f, B_M)~\leq~ \frac{\Lambda}{2}\,.
\end{equation}
 
In order to prove our claim, we decompose $\cal J$ in the countable union of suitable pieces. First of all, for every integer $k\geq 1$ we define
\[
{\cal J}_k~\doteq~ \left\{(\tau, \xi)\in {\cal J} : |u^+ (\tau, \xi)-u^- (\tau, \xi)| \geq \frac{1}{k}\right\}. 
\]
Obviously, ${\cal J} = \bigcup_k {\cal J}_k$. Next, given any pairs $(\tau_1, \xi_1)$ and $(\tau_2, \xi_2)$ in ${\cal J}_k$, consider the two piecewise constant functions
\begin{equation}\label{e:Uj}
U_j(t,x)~\doteq~
\left\{\bega{rl} u^- (\tau_j, \xi_j)\qquad &\hbox{if}\qquad x<\lambda (\tau_j, \xi_j) t,
\\[2mm] u^+(\tau_j,x_j) \qquad &\hbox{if}\qquad x>\lambda (\tau_j, \xi_j) t,\enda\right.
\qquad\qquad j=1,2.
\end{equation}
By \eqref{mspeed} there is a positive number $\varepsilon (k, \Lambda)$ depending on $k$ and $\Lambda$ such that the following holds. If $(\tau_s, \xi_s)$ is yet a third point in the plane with the properties that $\tau_s^2 + \xi_s^2 =1$ and $|\xi_s|\geq \Lambda |\tau_s|$, then  if we ``shift'' $U_2$ by this vector we get the inequality
\begin{equation}\label{e:lower-bound-2-states}
\int_{B_1} \bigl|U_1 (t,x) - U_2 (t+\tau_s, x+\xi_s)\bigr|\, dx\, dt ~\geq ~6 \varepsilon\, ,
\end{equation}
where $B_1$ denotes the unit disk centered at the origin in $\mathbb R^2$. We subdivide further ${\cal J}_k$ as a union of sets 
${\cal J}_{k,j}$, $j\geq 1$, where $(\tau, \xi)$ belongs to ${\cal J}_{k,j}$ if 
\begin{equation}\label{e:Jkj}
\int_{B_r (\tau, \xi)} |u (t,x) - U (t-\tau, x-\xi)|\, dx\, dt \leq \varepsilon r^2 \qquad \forall r<\frac{1}{j}\, .
\end{equation}
Here $U$ is defined as in (\ref{5}) and $B_r (\tau, \xi)$ denotes the disk centered at $(\tau, \xi)$ with radius~$r$. 
Clearly,
\[
{\cal J} ~=~ \bigcup_k {\cal J}_k ~=~ \bigcup_{k,j} {\cal J}_{k,j}\,.\]
Next, consider two points  $(\tau_1, \xi_1), (\tau_2, \xi_2)\in {\cal J}_{k,j}$ such that $r\doteq 
\bigl|(\tau_1, \xi_1) - (\tau_2, \xi_2)\bigr|< \frac{1}{2j}$.   Let $U_j$ be as in (\ref{e:Uj}) and  set the shift $(\tau_s, \xi_s)$ to be
$$
\tau_s \,=\, \frac{\tau_1-\tau_2}{r}\,,  \qquad\quad 
\xi_s \,=\, \frac{\xi_1-\xi_2}{r}\, .
$$
 We claim that \eqref{e:lower-bound-2-states} cannot hold with this shift.  This will enable us to conclude $|\xi_s|\leq \Lambda |\tau_s|$. 
To prove the claim, observe first that
\begin{align*}
& \int_{B_1} |U_1 (t,x) - U_2 (t+\tau_s, x+\xi_s)|\, dx\, dt\\ 
\leq &
\int_{B_1} |U_1 (t,x) - u (\tau_1 + rt, \xi_1 + rx)|\, dx\, dt + \int_{B_1} |u (\tau_1 + rt, \xi_1 + rx) - U_2 (t+\tau_s, x+\xi_s)|\, dx\, dt
\end{align*}
We then change variables in the integrals to $(\sigma, y) = (\tau_1 + rt, \xi_1+rx)$. Observe that
\[
\bega{rl}
t+\tau_s &=~ t + r^{-1} (\tau_1-\tau_2) ~=~ r^{-1} (t+\tau_1-\tau_2)= r^{-1} (\sigma - \tau_2),\\[1mm]
x+\xi_s &=~ x + r^{-1} (\xi_1-\xi_2) ~=~ r^{-1} (y-\xi_2) ,\enda
\]
while $U_j \bigl(r^{-1} (\sigma-\tau_j), r^{-1} (y-\xi_j)\bigr) = U_j (\sigma-\tau_j, y-\xi_j)$ because of the $0$-homogeneity of the functions $U_j$. Hence the change of variables yields
\[
\bega{l} \ds
 \int_{B_1} \bigl|U_1 (t,x) - U_2 (t+\tau_s, x+\xi_s)\bigr|\, dx\, dt\\[4mm]
\qquad\qquad \qquad \leq ~\ds\frac{1}{r^2} \int_{B_r (\tau_1, \xi_1)} \bigl|U_1 (\sigma-\tau_1, y-\xi_1) - u (\sigma,y)\bigr|\, dy\, d\sigma\\[4mm]
\qquad\qquad \qquad\qquad \qquad\qquad\qquad\ds+ \frac{1}{r^2} \int_{B_r (\tau_1, \xi_1)} \bigl|u (\sigma,y) - U_2 (\sigma -\tau_2, y-\xi_2)\bigr|\, dy\, d\sigma\\[4mm]
\ds \qquad\qquad\qquad\leq ~\frac{1}{r^2} \int_{B_r (\tau_1, \xi_1)} \bigl|U_1 (\sigma-\tau_1, y-\xi_1) - u (\sigma,y)\bigr|\, dy\, d\sigma
\\[4mm]
\qquad\qquad\qquad\qquad\qquad\qquad\qquad\ds + \frac{1}{r^2} \int_{ B_{2r} (\tau_2, \xi_2)}
\bigl|u (\sigma ,y) - U_2 (\sigma-\tau_2, y-\xi_2)\bigr|\, dy\, d\sigma,
\enda
\]
where we have used the inclusion $B_r (\tau_1, \xi_1)\subset B_{2r} (\tau_2, \xi_2)$ to get the last inequality. Note next that $r<2r< \frac{1}{j}$ and, since both $(\tau_1, \xi_1)$ and $(\tau_2, \xi_2)$ belong to ${\cal J}_{k,j}$, we can use (\ref{e:Jkj}) to bound the first summand by $\varepsilon$ and the second summand by $4\varepsilon$. In particular we conclude
\[
\int_{B_1} \bigl|U_1 (t,x) - U_2 (t+\tau_s, x+\xi_s)\bigr|\, dx\, dt \leq 5 \varepsilon\, .
\]
As already pointed out, since the latter inequality contradicts (\ref{e:lower-bound-2-states}), we conclude that $|\xi_s|\leq \Lambda |\tau_s|$, which in turn gives $|\xi_2-\xi_1|\leq \Lambda |\tau_2-\tau_1|$.

We have thus proved the following fact:
\begin{itemize}
\item[{\bf (L)}] 
{\it If $(\tau_1 , \xi_1), (\tau_2 , \xi_2)\in {\cal J}_{k,j}$ and $|(\tau_1, \xi_1) - (\tau_2, \xi_2)|< \frac{1}{2j}$, then $|\xi_2-\xi_1|\leq \Lambda |\tau_2-\tau_1|$.}
\end{itemize} 
It is well known that from {\bf (L)} it follows that ${\cal J}_{k,j}$ can be covered by countably many Lipschitz graphs of functions $x=\phi_\ell (t)$, $\ell\in \mathbb N$, with Lipschitz constant at most $\Lambda$.  See for instance \cite[Lemma 4.7]{DL}. 

For readers' convenience, we include here a proof. If $B = B_{1/{4j}} (x_0, t_0)$ is any disk of radius $1/(4j)$ and we set $F\doteq {\cal J}_{k,j}\cap B$, then {\bf (L)} implies
\begin{equation}\label{e:Lip}
|x-y|\leq \Lambda |t-s| \qquad \forall (t,x), (s,y)\in F\, .
\end{equation}
This obviously implies that there are no points of $F$ which lie on the same line $\{t= {\rm const}\}$. Hence, if $G$ is the projection of $F$ on the time axis, then there is a function $\phi: G \to \mathbb R$ such that $F = \{(t, \phi (t)): t\in G\}$. On the other hand \eqref{e:Lip} is equivalent to the statement that the Lipschitz constant of $G$ is at most $\Lambda$. By the classical Lipschitz extension theorem we can simply extend $\phi$ to a Lipschitz function defined on the whole time axis. Since ${\cal J}_{k,j}$ can be covered by a countable collection of disks with radius $1/(4j)$, the existence of the desired covering by means of countably many Lipschitz graphs
follows immediately. 
\v
{\bf 2.} Next, we wish to show that, if $t\not\in {\cal N}$ and $(t, x)\not \in \mathbf{J}$, then $y\mapsto u (t, y)$ is continuous at $x$. 
We start by noticing that, since $(t,x)\not \in \mathbf{J}$, $u$ is approximately continuous at $(t,x)$ as a function of two variables. Therefore there is a $u_0\in \mathbb R^n$ such that 
\[
\lim_{r\downarrow 0} \frac{1}{r^2} \int_{t-r}^{t+r} \int_{x-r}^{x+r} |u (s,y) - u_0|\, ds\, dy ~= ~0\, . 
\]
In particular, for every fixed $\varepsilon > 0$ there is an $r_0 (\varepsilon) >0$ such that 
\[
 \int_{t-r}^{t+r} \int_{x-r}^{x+r} |u (s,y) - u_0|\, ds\, dy ~\leq ~\varepsilon r^2 \qquad \forall r \leq r_0 (\varepsilon)\, . 
\]
An elementary application of Chebyshev's inequality and Fubini's theorem yields then the existence of a $t (r)\in [t - r\sqrt{\varepsilon}, t+r \sqrt{\varepsilon}]$ such that 
\[
\int_{x-r}^{x+r} |u (t (r), y) - u_0|\, dy ~\leq ~r \sqrt{\varepsilon}\, .
\]
Furthermore we can use the Lipschitz estimate (\ref{ulip}) to bound
\begin{equation}\label{e:1d-approx-continuity}
\frac{1}{r} \int_{x-r}^{x+r} |u (t, y) - u_0|\, dy \leq \sqrt{\varepsilon} + C_M \frac{|t-t(r)|}{r} 
~\leq ~(1+C_M) \sqrt{\varepsilon} \qquad \forall r \leq r_0 (\varepsilon)\, .
\end{equation}
On the other hand recall that $y\mapsto u (t, y)$ is a function of bounded variation on the real line. As such, every point is either a classical jump point, or a point of continuity. Since $\varepsilon$ in \eqref{e:1d-approx-continuity} can be closen arbitrarily small,  $x$ cannot be a classical jump point of $u (t, \cdot)$ and must therefore be a point of continuity, which was in fact our initial claim.
\v
{\bf 3.} Together with the functions $\phi_\ell$ in (\ref{p1}),
we consider functions of the form
\bel{p2} \phi^{y+}(t) = y+\Lambda t,\qquad\qquad \phi^{y-}(t) = y-\Lambda t,\qquad\qquad y\in {\mathbb Q}.\eeq
Since here $y$ is a rational point, there are countably many of these functions.
For convenience, the countable set of all functions in
(\ref{p1}) together with those in (\ref{p2}) will be relabeled as
\bel{sin}
\{ \psi_n\,;~~n\geq 1\}.\eeq 

Next, we observe that, for every  $j,k\geq 1$, the 
scalar function
\bel{Wpm} W_{jk} (t)~\doteq~\left\{ \bega{cl} \tv\big\{ u(t)\,;~]\psi_j(t)\,,~\psi_k(t)[\,\bigr\}
\qquad &\hbox{if}\quad \psi_j(t)<\psi_k(t)\,,\\[2mm]
0\qquad &\hbox{otherwise,}\enda\right.\eeq
is bounded and measurable (indeed, it is lower semicontinuous).  Therefore a.e.~$t\in [0,T]$ is a Lebesgue point.
We denote by ${\cal N}'\subset [0,T]$ the set of all times $t$ which are NOT
Lebesgue for at least one of the countably many functions  $W_{jk}$.
Of course, ${\cal N}'$ has zero Lebesgue measure.

In view of (\ref{errest}), we will prove the theorem by establishing  the following claim.
\begi
\item[{\bf (C)}] {\it For every $\tau\in [0,T]\setminus ({\cal N\cup \cal N'})$ and $\ve>0$,
one has}
\bel{limsup}
\limsup_{h\to 0+} {1\over h} \Big\| u( \tau + h) - S_h u( \tau)\Big\|_{\L^1}~\leq~\ve.\eeq
\endi 
 \v
 {\bf 4.} Assume $\tau \notin \cal N\cup \cal N'$.  
By induction on $k=0,1,2,\ldots, N$,
we will construct points
$$y_0\leq y_0'~<~y_1\leq y'_1~<~y_2\leq y'_2~<~\cdots~<~y_N\leq y_N'$$
with the following properties.
\begi
\item[(i)]  Either $y_k = y_k' = \phi_j(\tau)$ for some $j$, or else $y_k<y'_k$ and
\bel{rational}y_k'+\Lambda \tau\in {\mathbb Q},\qquad\qquad y_k-\Lambda \tau\in {\mathbb Q}\,.\eeq
\item[(ii)]  The first and the last points satisfy
\bel{tve1}
\tv\bigl\{ u( \tau,\cdot)\,;~]-\infty, \,y_0'[\,\bigr\}~<~\ve,
\eeq
\bel{tve11} \tv\bigl\{ u( \tau,\cdot)\,;~]y_N,\, +\infty[\,\bigr\}~<~\ve.
\eeq
Moreover, 
for every $k\in\{1,\ldots,N\}$,  considering the total variation of the   right-continuous function $u(\tau,\cdot)$
on the following open and half-open intervals, one has
\bel{tv3} 
\tv\bigl\{ u(\tau, \cdot)\,;~]y_{k-1}, \,y_k'[\,\,\bigr\}~<~2\ve,
\eeq
\bel{tv4} 
\tv\bigl\{ u(\tau, \cdot)\,;~]y_{k-1}, \,y_k']\,\bigr\}~>~\ve,
\eeq
\bel{tv5} 
\tv\bigl\{ u(\tau, \cdot)\,;~]y_{k}, \,y_k']\,\,\bigr\}~<~{\ve\over 2}\,.
\eeq

\endi
\v
The construction is straightforward.  We first determine  points $y_0<y_0'$ 
satisfying (\ref{rational}) and  such that (\ref{tve1}) holds together with 
$$ \tv\bigl\{ u( \tau,\cdot)\,;~]y_0\,,\, y_0']\,\bigr\}~<~{\ve\over 2}\,.$$

Next, assume by induction that the points $y_0,\ldots y_{k-1}, y_{k-1}'$ have
already  been constructed.  Consider the point
$$z_{k} ~=~\sup~\Big\{ x>y_{k-1}'\,;~~\tv\bigl\{ u(\tau,\cdot)\,;~~]y_{k-1}, x[\,\bigr\} <{3\ve\over 2}\Big\}.
$$
CASE 1:   If the map $y\mapsto u(\tau,y)$ has a jump at $z_k$, then by step {\bf 2} we have
$(\tau, z_k)\in {\bf J}$.    Hence by step {\bf 1} it follows
$z_k = \phi_j(\tau)$ for some $j$.
In this case we set $y_k=y_k'=z_k$.
\v
CASE 2:  If the map $y\mapsto u(\tau,y)$ is continuous at $z_k$, then we can take two points
$y_k<z_k<y_k'$ such that (\ref{rational}) holds,  together with (\ref{tv3})--(\ref{tv5}).
\v
Since $\tv \bigl\{ u(\tau,\cdot)\bigr\}\leq M$,  by (\ref{tv4})
the total number of
points $y_k$ will be 
\bel{NM} N~\leq~{2M\over \ve}\,.\eeq

\begin{figure}[htbp]
\centering
 \includegraphics[scale=0.45]{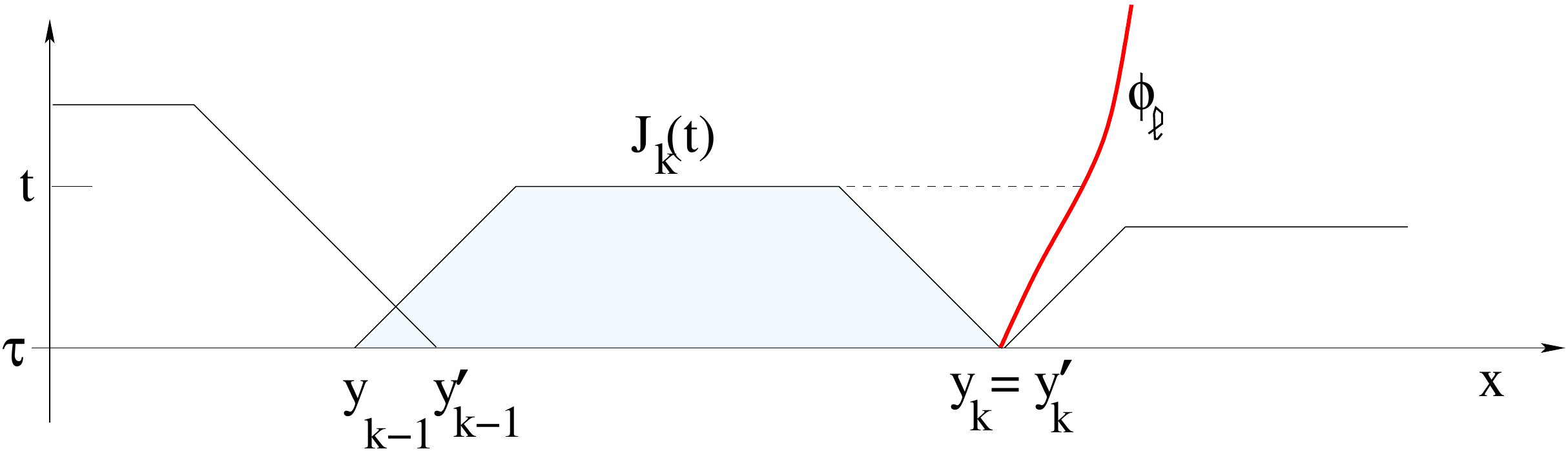}
    \caption{\small  The points $y_k$, $y_k'$ constructed 
 in the proof of the theorem.  
}
\label{f:hyp252}
\end{figure}

\v
\n{\bf 5.}  The remainder of the proof is very similar to the one in \cite{BGu}.
For any given $y\in \R$, we denote by $U^\sharp= U^\sharp_{(u,\tau,y)}(t,x)$ the solution to the 
Riemann problem for (\ref{1}), with initial data
\bel{Rdata}u(\tau, x)~=~\left\{ \bega{rl} u(\tau, y-)\qquad &\hbox{if} ~~x< y,\cr
u(\tau, y+)\qquad &\hbox{if} ~~x> y.\enda\right.\eeq
Moreover, for every given $k\geq 1$ we denote by $U^\flat= U^\flat_{(u,\tau,k)}(t,x)$ the solution to the 
linear Cauchy problem with constant coefficients
\bel{linh} v_t + A v_x~=~0,\qquad \qquad  v(\tau,x) = u(\tau, x).\eeq
Here the $n\times n$ matrix $A$ is the Jacobian matrix of $f$ computed
at the midpoint of the interval $[y_{k-1}, y'_k]$. Namely, 
$$A~=~ Df\left(u\Big(\tau, {y_{k-1}+ y'_k\over 2}\Big)\right).$$
 With reference to Fig.~\ref{f:hyp252}, to estimate  the lim-sup
in (\ref{limsup}), we need to estimate two types of integrals.

\begi
\item[(I)] The integral of $\bigl| u(t,x) - U^\sharp_{(u,\tau, y_k)} (t,x)\bigr|$ over the intervals
$$J_y(t)~\doteq~
\bigl[  y_k -\Lambda (t-\tau)\,,~ y_k +\Lambda(t-\tau) \bigr],$$
for every $k$ such that $y_k= y_k'$.

\item[(II)]The integral of $\bigl| u(t,x) - U^\flat_{(u,\tau,k)} (t,x)\bigr|$  over the intervals
$$J_{k}(t)~=~\bigl[ y_{k-1}+\Lambda (t-\tau)\,,~
 y_{k}'- \Lambda(t -\tau) \bigr].$$

\endi
%
%
%
%
%
\v
{\bf 6.}  To estimate integrals of type (I) we observe that, since $\tau\notin \N$, 
$(\tau, y_k)$ is either
a Lebesgue point or  a point of approximate jump of the 
function $u$. Therefore
\bel{Ush}\lim_{t\to \tau+} ~{1\over t-\tau}\int_{y_k-\Lambda (t-\tau)}
^{y_k+\Lambda (t-\tau)}\Big| u(t,x)- U^\sharp
_{(u;\tau,y_k)}(t, x)\Big|~dx~=~0.\eeq
Indeed, this follows from (\ref{6}) and the Lipschitz continuity of the map $t\mapsto u(t,\cdot)\in \L^1(\R)$.   See Theorem~2.6 in \cite{Bbook} for details. 
\v
{\bf 7.} To estimate integrals of type (II),   
two main cases will be considered.

CASE 1: $y_{k-1}-\Lambda \tau\in{\mathbb Q}$ and $y'_{k}+\Lambda \tau\in{\mathbb Q}$.
In this case, since $\tau\notin {\cal N}'$, the function
\bel{Vtdef}V(t)~\doteq~\tv\Big\{ u(t,\cdot)\,;~\bigl]y_{k-1}+\Lambda (t-\tau)\,,~y'_{k}-\Lambda (t-\tau)\bigr[\Big\}\eeq
has a Lebesgue point at $t=\tau$. Hence
$$\lim_{h\to 0+} {1\over h} \int_\tau^{\tau+h} \bigl| V(t)-V(\tau)\bigr|\, dt~=~0.$$
Since $V(\tau)\leq 2\ve$, this implies
\bel{lsuV}\limsup_{h\to 0+} {1\over h} \int_\tau^{\tau+h} V(t)\, dt~\leq~2\ve.\eeq
\v
CASE 2:  $y_{k-1} = \phi_\ell(\tau)$, $y'_k = \phi_{\ell'}(\tau)$ 
for some indices $\ell, \ell'\in {\mathbb N}$. 
In this case, since $\tau\notin {\cal N}'$,  the function 
$$W(t)~\doteq~\tv\Big\{ u(t,\cdot)\,;~\bigl]\phi_\ell(t)\,,~\phi_{\ell'}(t)\bigr[\Big\}$$
has a Lebesgue point at $t=\tau$.
Recalling that the functions $\phi_\ell$, $\phi_{\ell'}$ have Lipschitz constant $\leq \Lambda$, 
a comparison with (\ref{Vtdef}) immediately yields $V(t)\leq W(t)$ for all $t\geq \tau$.  
Since our construction implies
$W(\tau)\leq 2\ve$, we thus conclude that (\ref{lsuV}) again holds.

The remaining two cases, where $y_{k-1} = \phi_\ell(\tau)$ and  $y'_k+\Lambda \tau\in{\mathbb Q}$,   or where $y_{k-1}-\Lambda \tau\in{\mathbb Q}$ and $y'_k = \phi_{\ell'}(\tau)$ for some $\ell,\ell'\in
{\mathbb N}$,
can be handled in the same way.  Namely, (\ref{lsuV}) always holds.
\v
Using again the fact that the map 
$t\mapsto u(t,\cdot)$ is Lipschitz continuous from $[0,T]$ into $\L^1(\R)$,
the same argument used in \cite{BGu} now yields
\bel{Ubes}
  \limsup_{t\to\tau+}
  \frac{1}{t-\tau}  \int_{y_{k-1} + \Lambda (t-\tau)}^{y_k'-\Lambda (t-\tau)} \left|u(t,x)-U^{\flat}_{(u,\tau,k)}\left(t,x\right)\right|\;
dx~=~\O(1)\cdot\varepsilon^{2}.
\eeq
Indeed, this corresponds to formula (3.20) in \cite{BGu}.  Based on (\ref{lsuV}), the proof is identical and will not be repeated here.
\v
{\bf 8.} On the other hand, as showed in \cite{BiB} all trajectories of the semigroup are weak solutions of (\ref{1}) which satisfy the Liu admissibility conditions.  Therefore they
satisfy the same bounds as in (\ref{Ush}) and (\ref{Ubes}).   More precisely:
\bel{S1}\lim_{t\to \tau+} {1\over t-\tau}\int_{y_k-\Lambda (t-\tau)}
^{y_k+\Lambda (t-\tau)}\Big| \bigl(S_{t-\tau} u(\tau)\bigr) (x)- U^\sharp_{(u;\tau,y_k)}(t, x)\Big|~dx~=~0,\eeq
\bel{S2} \limsup_{t\to\tau+}
  \frac{1}{t-\tau}  \int_{y_{k-1} + \Lambda (t-\tau)}^{y_k'-\Lambda (t-\tau)} \Big|  
  \bigl(S_{t-\tau} u(\tau)\bigr) (x)- U^\flat_{(u,\tau,k)}(t, x)\Big|~dx~=~\O(1)\cdot \ve^2.\eeq
\v
{\bf 9.}  Combining  the previous estimates (\ref{Ush}), (\ref{Ubes}), (\ref{S1}), (\ref{S2}), and recalling that the total number of 
intervals is $N\leq  2M\ve^{-1}$, we establish the limit (\ref{limsup}), proving the theorem.
\endproof

\v

{\bf Acknowledgments.} The research by the first author
 was partially supported by NSF with
grant  DMS-2006884, ``Singularities and error bounds for hyperbolic equations".

\end{document}